\documentclass[11pt,a4paper]{article}
\usepackage[margin=2 cm]{geometry}
\usepackage{float}
\usepackage{times}
\usepackage{graphicx}
\usepackage{amssymb}
\usepackage{mathtools}
\usepackage{enumitem}
\usepackage{gensymb}
\usepackage{subcaption} 
\usepackage{array}
\makeatletter
\let\NAT@parse\undefined
\makeatother
\usepackage[backend=biber,style=numeric-comp, sorting=none, isbn=false, doi=true, url=false, giveninits=true]{biblatex}
\usepackage{makecell}
\usepackage{longtable}
\usepackage{amsmath, mathtools, amsfonts}
\usepackage{amsthm}
\usepackage{tablefootnote}
\usepackage{booktabs}
\usepackage{tikz}
\usepackage{lipsum}
\usepackage{url}
\usetikzlibrary{positioning,fit,backgrounds, calc}
\usetikzlibrary{shapes.geometric, arrows}
\usepackage{pgfplotstable} 
\pgfplotsset{compat=1.18}
\usepackage{xstring}

\usepackage{todonotes}

\usepackage[colorlinks=true, allcolors=blue]{hyperref} 
\usepackage{orcidlink}

\usepackage{authblk}
\author[1]{Janina Zittel \orcidlink{0000-0002-0731-0314}}
\author[2]{Annika Buchholz \orcidlink{0009-0008-3235-0896}}
\author[3]{Michael Bussieck}
\author[4]{Frederik Fiand}
\author[5]{Thorsten Koch \orcidlink{0000-0002-1967-0077}}
\author[6]{Lukas Mehl}
\author[7]{Niels Lindner  \orcidlink{0000-0002-8337-4387}}
\author[8]{Manuel Wetzel}

\affil[1,2,5,6]{Zuse Institute Berlin, Berlin, Germany, \texttt{zittel@zib.de (CA), buchholz@zib.de, koch@zib.de, mehl@zib.de}}
\affil[3,4]{GAMS Software GmbH, Frechen, Germany, \texttt{mbussieck@gams.com, ffiand@gams.com}}
\affil[7]{Freie Universität Berlin, Germany, \texttt{lindner@zib.de}}
\affil[8]{Deutsches Zentrum für Luft- und Raumfahrt e.V., Stuttgart, Germany, \texttt{manuel.wetzel@dlr.de}}

\title{Computational acceleration strategies for large-scale energy system optimization: \\
a comparative study of GPU-accelerated \\
and distributed-memory solvers}

\begin{document}
\maketitle

\abstract{\small Energy system optimization models are increasing in scope and resolution, integrating detailed technology representations, sector coupling, and multiple scenarios reflecting uncertainty in future energy demand. These advances yield large and challenging linear programs whose efficient solution remains a bottleneck in practical energy system analyses. For a long time, the standard way to address such problems has relied on shared-memory interior-point methods (IPM), which combine robustness and accuracy but face scalability limits as model instance size grows.

Specialized solver architectures are beginning to change this picture. Two promising directions have emerged: (i) GPU-accelerated first-order methods (FOM) such as the primal-dual linear programming approach; and (ii) distributed-memory IPM, exemplified by the open-source solver PIPS-IPM++, which can exploit block structure that arises in many energy system models and achieves high parallelism on high-performance computing systems. These developments open new opportunities for large-scale optimization in energy system analysis.

This paper presents a computational study comparing these solver classes on a diverse test set of large-scale linear programs arising from energy system analysis, including scenario-based formulations derived from stochastic programming. We investigate how parallelization strategies affect solution time and numerical accuracy. The results illustrate that distributed-memory IPM can leverage problem structure to deliver substantial speed-ups on specific problems with block-angular structures. GPU-accelerated FOMs demonstrate strong scalability but may yield solutions with higher relative infeasibilities, which, depending on the use case and model uncertainty, can still be acceptable.

Overall, our findings indicate that recent algorithmic and hardware advances substantially broaden the computational toolbox available to the energy system optimization community. Each solver class exhibits distinct advantages: shared-memory IPMs remain a powerful tool for reliably obtaining high-accuracy solutions; distributed-memory IPMs can extend scalability to hundreds of cores for certain structured models, enabling faster time-to-solution; and GPU-based FOM can deliver fast solutions when such lower accuracy levels are appropriate. Together, they help make high-resolution, multi-scenario energy system optimization models tractable across a broader range of problem sizes and computing environments.}

\section{Introduction} \label{Introduction} 
Energy supply systems are undergoing profound structural and regulatory transformations, shifting toward highly decentralized, multi-sectoral, and interconnected networks that integrate an increasing share of variable renewable energy sources. These changes are reflected in energy system optimization models (ESOM), which are increasingly high-dimensional and interconnected, capturing the complexity of continental-scale, multi-sectoral energy transition pathways \parencite{DeCarolis2017}. Most optimization models are still formulated under the assumption of a deterministic system with perfect foresight. However, robust decision-making increasingly requires accounting for uncertainty, including the stochastic nature of renewable generation, weather-dependent demand fluctuations, and uncertain technology costs. This requires solving these massive systems repeatedly across numerous scenarios or other approaches such as stochastic modelling \parencite{YUE2018204, Loeffler2022,  Lindneretal2025, Muschner2024}. As these models integrate detailed technology representations and expansive stochastic scenarios to account for weather-driven uncertainty and cost fluctuations, the resulting linear programs (LP) often exceed hundreds of millions of variables. As a result, model instance sizes increasingly exceed the memory limits of traditional single-node computing architectures. 

For decades, the field has relied on mathematical solvers that have seen significant improvements \parencite{Koch2022}, with shared-memory interior-point methods (IPM), which typically implement a barrier algorithm, now representing the fastest option for most large-scale LP formulations. However, the linear algebra operations at the core of IPMs require large amounts of memory and scale poorly across distributed computing systems. 

To mitigate this computational burden, researchers typically employ one of two different strategies: model-based acceleration or solver-based speed-up methods \parencite{Cao2019, Kotzur2021}. Model-based approaches often utilize temporal or spatial aggregation to simplify the problem, or model-centric mathematical decomposition methods such as Benders decomposition, stochastic dual dynamic programming, or Lagrangian relaxation. However, temporal aggregation techniques are unsuitable for systems characterized by high shares of decentralized renewables and long-duration storage \parencite{Kotzur2018},  which are precisely the most relevant system configurations for modern decarbonization studies. While model decomposition techniques such as Benders decomposition, column generation, or Lagrangian relaxation can yield substantial computational gains, they introduce their own challenges: convergence can be sensitive to problem structure, implementation requires significant algorithmic expertise and often extensive remodeling efforts, and the resulting solutions may still face scalability limits for the most ambitious studies. Moreover, for very complex analyses -- such as continental-scale, sector-coupled systems with high temporal resolution and extensive uncertainty quantification -- even sophisticated model reduction approaches must often be combined with solver-level acceleration strategies to overcome computational boundaries. 

Fortunately, recent advances in optimization technology are beginning to expand our repertoire of methods through two complementary solver-based pathways. First, distributed-memory IPMs extend the reach of traditional IPMs by exploiting the inherent block-diagonal structure that naturally arises in many ESOM formulations, particularly in multi-scenario and multi-stage stochastic programs. The open-source solver PIPS-IPM++ \parencite{Rehfeldtetal22, Kempke2025} exemplifies this approach, distributing computational workload across HPC cluster nodes to achieve high parallelism while maintaining the numerical reliability and precision of IPMs. Although a systematic assessment of which ESOM classes benefit most from distributed-memory parallelization remains an open research question, preliminary studies with the ESOMs REMix \parencite{WETZEL2025101893} and oemof-B3 \parencite{KochURB2025} have demonstrated promising scalability for certain problem structures. Second, GPU-accelerated first-order methods (FOM) leverage the massive parallelism of modern GPUs to overcome the memory and scalability limitations of IPMs. Among these, the primal-dual linear programming (PDLP) approach has seen substantial algorithmic refinement in recent years \parencite{Applegate25, HPR-LP2025, ZhangBoyd25}, with implementations such as NVIDIA's open-source solver cuOpt and variants, such as the Primal-Dual Hybrid Gradient method, increasingly integrated into commercial and academic solver suites. FOMs have demonstrated promise not only for large-scale LP but also as components of advanced mixed-integer programming (MIP) solution strategies \parencite{KempkeKoch25}. Together, these architectures provide new opportunities to solve next-generation ESOMs without reducing spatial or temporal resolution. 

This study evaluates how emerging solver technologies perform on realistic large-scale ESOMs. To this end, we provide a computational study benchmarking three distinct solver classes: traditional shared-memory IPMs, distributed-memory IPMs, and GPU-accelerated FOMs, including three open-source solvers: HiGHS\cite{Qi2018},  PIPS-IPM++\cite{Rehfeldtetal22} and  NVIDIA cuOpt\cite{cuopt}. Our test suite comprises a diverse set of large-scale LPs derived from real-world energy system analysis. We investigate how various parallelization strategies impact critical performance metrics, including solution time and numerical accuracy. While our primary focus remains on LP formulations, these findings hold significant implications for MIP. In cases where MIPs are ``LP-expensive'' -- i.e., the relaxation consumes the bulk of the runtime -- our results can be integrated into hybrid strategies \parencite{KempkeKoch25} and MIP frameworks such as the Ubiquity Generator, designed to accelerate the integer search via optimized LP sub-problems \parencite{UG}.

A second objective of this work is to contribute to transparent benchmarking within the energy research community. Consequently, all instances used in this study that permit public sharing are consolidated into a Zenodo database to facilitate future benchmarking activities. This effort complements existing initiatives in the energy community, such as the Open Energy Transition benchmarking suite, by specifically extending the problem size and complexity to reflect the challenges of the next generation of models. Ultimately, this experimental study provides actionable insights for modelers to determine which solution strategies -- distributed architectures or FOM -- are worth exploring for specific problem topologies, rather than defaulting to the traditional combination of shared-memory IPMs and model complexity reduction.

\section{Modeling Base} \label{section:problem}
We evaluate the performance of our selected acceleration strategies using a benchmark suite of 46 large-scale LP instances. These problems represent the computational frontier of modern energy system analysis, where high spatial and temporal resolutions often lead to dimensions that challenge the limits of traditional shared-memory architectures.

To ensure the benchmark is representative of modern modeling practice, the test suite distinguishes between two primary functional categories: operational dispatch and investment-oriented capacity expansion. In dispatch models, the focus lies on the cost-optimal scheduling of existing assets, where inter-temporal coupling is primarily driven by energy storage constraints, such as state-of-charge limits. In contrast, investment models include endogenous decisions for asset sizing, which introduces a much more demanding coupling structure. Because a single capacity variable appears in the constraints of every operational hour, investment decisions inherently couple all time steps. This global coupling significantly increases the complexity of the constraint matrix and serves as a primary driver for the emergence of dense columns, which pose a specific challenge for IPMs.

The suite covers the two main categories in energy system optimization. On one hand, more traditional integrated energy system models such as TIMES and SWITCH, and on the other hand, models derived from the power sector modeling with hourly resolution, such as PyPSA, REMix, ETHOS.FINE, genx. The main difference is typically the temporal resolution. Instances such as JRC-EU-TIMES-dispatch-30-2016ts, ethos fine europe 60tp-175-720ts, and SWITCH-China-open-model 32-433ts represent continental-scale systems with up to 175 regional nodes (ETHOS.FINE) on up to 2016  time slices (JRC-EU-TIMES). The suite also includes two-stage stochastic models, such as the oemofB3\_int series, which encompass up to 250 scenarios in a single LP. Crucially, all instances in this study are native LP formulations. Unlike benchmarks that rely on the root-node relaxations of MIP, these models are designed as pure LPs to support large-scale planning where integrality is often computationally prohibitive.

To understand the problem structure, we analyzed dimensions as well as row and column densities as proxies for constraint and variable coupling. As reported in Table \ref{tab:instances}, instance dimensions range from $10^6$ to over $2 \times 10^8$ rows and columns, with non-zeros (NZ) spanning 27 million to 2.68 billion in the BEAM\_4032\_11\_8\_CLI instance. These metrics directly impact decomposability: high NZ counts in columns indicate variables coupling numerous constraints, while dense rows represent constraints linking many variables. Consequently, strong coupling makes even smaller instances computationally challenging. Our test set reflects this spectrum: zen-garden and oemof are highly sparse, with a median of 2 NZ per row/column. Conversely, industrial\_[2-6] exhibit extreme coupling, with the densest 2~\% of rows and columns containing over 200 and 2,000 NZ, respectively.


The test suite consists of 46 instances categorized by their accessibility. Ten previously unpublished instances from various frameworks have been consolidated into \textit{Zenodo Project 18953379} \cite{zenodo_neu} to facilitate future benchmarking activities. 22 of \textit{reference instances} provide a baseline for comparison. Among others, these have been published by the Open Energy Transition Benchmark \cite{OET}. 14 \textit{confidential industrial instances} are derived from real-world grid planning data. While proprietary regulations prohibit their public release, they are included to validate performance on data featuring different numerical scaling and coupling patterns than research models.

\begin{table}[htbp]
    \footnotesize
    \centering
        \caption{Overview on the ESOM instances, models, and linear program properties}
    \begin{tabular}{ll rrr}
   
        \hline
         \textbf{short} & \textbf{model } & \textbf{rows} & \textbf{columns} & \textbf{non-zeros} \\
         \hline
         
         \multicolumn{5}{l}{\textbf{zenodo.18953379 \cite{zenodo_neu}:}} \\
        288\_22\_8 & TIMES & 16{,}247{,}539     & 11{,}833{,}217     & 93{,}254{,}832   \\  
        672\_22\_8\_cli & TIMES & 37{,}645{,}233     & 27{,}330{,}735     & 215{,}840{,}892 \\
        BEAM\_2016\_11\_8 & TIMES & 63{,}287{,}007     & 45{,}315{,}135     & 369{,}638{,}727 \\   
        BEAM\_4032\_11\_8\_CLI             & TIMES & 171{,}349{,}731    & 116{,}118{,}348    & 2{,}673{,}203{,}562 \\
        ELMOD\_876\_10\_noVEnames & ELMOD & 256{,}284{,}723    & 226{,}061{,}766    & 717{,}436{,}984 \\
        nt2030\_2030\_CY2009\_st1 & openTEPES & 16{,}867{,}670     & 19{,}742{,}599     & 49{,}466{,}634  \\
        oemofB3\_int\_125 & oemof & 107{,}310{,}059 & 133{,}590{,}100    & 328{,}437{,}616 \\   
        oemofB3\_int\_250 & oemof & 214{,}620{,}059    & 267{,}180{,}100 & 656{,}875{,}116 \\    
        oemofB3\_int\_37 &  oemof & 31{,}763{,}819     & 39{,}542{,}740     & 97{,}217{,}616 \\ 
        OTAI\_b-ts5\_cplex & TIMES & 26{,}394{,}169     & 43{,}518{,}661     & 164{,}834{,}985 \\
        \hline
        \multicolumn{5}{l}{\textbf{reference instances:}} \\
         spineopt-multiyear-invest-eu-case-study$^{\mathrm{a}}$  & SpineOpt & 18{,}422{,}757     & 13{,}134{,}385     & 41{,}956{,}094 \\    
        remix\_unseen\_adalpert\_1h$^{\mathrm{b}}$ & REMix & 91{,}226{,}059     & 94{,}633{,}107     & 293{,}062{,}856 \\
        remix\_nagsys\_eu\_ineq\_1h$^{\mathrm{b}}$ & REMix & 215{,}009{,}179    & 199{,}532{,}154    & 877{,}558{,}092  \\
        remix\_unseen\_gunilpert\_1h$^{\mathrm{b}}$ & REMix &  38{,}857{,}052     & 40{,}635{,}185     & 130{,}896{,}035  \\
        remix\_nagsys\_cwe\_3h$^{\mathrm{b}}$ & REMix & 23{,}056{,}217     & 23{,}849{,}054     & 102{,}597{,}816\\
        remix\_nagsys\_eu\_8h$^{\mathrm{b}}$ & REMix & 24{,}817{,}865     & 25{,}472{,}312     & 108{,}677{,}973     \\ 
        remix\_nagsys\_fr\_1h$^{\mathrm{b}}$ & REMix & 16{,}102{,}023     & 16{,}651{,}889     & 69{,}463{,}361 \\
        remix\_nagsys\_cwe\_ineq\_1h$^{\mathrm{b}}$ & REMix & 76{,}614{,}493     & 71{,}163{,}989     & 314{,}482{,}461 \\
        remix\_yssp\_disp\_488r\_1h$^{\mathrm{b}}$ & REMix & 60{,}665{,}309     & 74{,}465{,}253     & 185{,}938{,}670 \\        
        ethos\_fine\_europe\_60tp-175-720ts$^{\mathrm{c}}$ & ETHOS.FINE & 9{,}230{,}938 & 8{,}141{,}182 & 32{,}766{,}436 \\
        genx-elec\_co2-15-168h$^{\mathrm{c}}$ & GenX & 11{,}832{,}633 & 13{,}777{,}202 & 42{,}847{,}400 \\
        genx-elec\_trex-15-168h$^{\mathrm{c}}$ & GenX & 11{,}832{,}658 & 13{,}777{,}227 & 41{,}476{,}545 \\
        genx-elec\_trex\_co2-15-168h$^{\mathrm{c}}$ & GenX & 11{,}832{,}689 & 13{,}777{,}228 & 43{,}263{,}428 \\
        JRC-EU-TIMES-dispatch-30-2016ts$^{\mathrm{c}}$ & TIMES & 19{,}401{,}352 & 23{,}208{,}969 & 366{,}217{,}330 \\
        pypsa-de-elec-50-1h$^{\mathrm{c}}$ & PyPSA & 18{,}370{,}417     & 8{,}637{,}885      & 35{,}860{,}746 \\    
        pypsa-de-sec-20-1h$^{\mathrm{c}}$    & PyPSA & 45{,}553{,}984     & 21{,}472{,}540     & 106{,}717{,}276  \\  
        pypsa-eur-elec-100-3h$^{\mathrm{c}}$  & PyPSA & 14{,}064{,}375     & 6{,}720{,}092      & 27{,}072{,}974  \\
        SWITCH-China-open-model\_32-433ts$^{\mathrm{c}}$ & SWITCH & 14{,}588{,}407 & 13{,}338{,}520 & 46{,}203{,}133\\
        TIMES-GEO-global-base-31-20ts$^{\mathrm{c}}$ & TIMES& 6{,}027{,}253 & 6{,}449{,}851 & 36{,}624{,}474 \\
        TIMES-GEO-global-netzero-31-20ts$^{\mathrm{c}}$ & TIMES & 5{,}063{,}755 & 5{,}422{,}595 &  30{,}487{,}085 \\
        zen-garden-eur-PI-ann-emis-lim-28-100ts$^{\mathrm{c}}$ & ZEN-garden & 16{,}737{,}720 & 14{,}342{,}280 & 47{,}988{,}391 \\
        zen-garden-eur-PI-constr-exp-28-100ts$^{\mathrm{c}}$ & ZEN-garden & 16{,}738{,}346 & 14{,}342{,}280 & 48{,}304{,}103 \\
        \hline
        \multicolumn{5}{l}{\textbf{confidential industrial instances}:} \\
        industrial\_01 & & 16{,}841{,}640     & 24{,}590{,}574     & 67{,}365{,}955\\ 
        industrial\_02 & & 1{,}045{,}745      & 1{,}874{,}896      & 148{,}620{,}736 \\
        industrial\_03 & & 1{,}050{,}757      & 1{,}886{,}886      & 150{,}019{,}065\\
        industrial\_04 & & 1{,}053{,}309      & 18{,}999{,}31      & 149{,}894{,}072 \\
        industrial\_05 & & 1{,}026{,}568      & 1{,}830{,}777      & 156{,}497{,}994 \\
        industrial\_06 & & 1{,}019{,}855      & 1{,}827{,}186      & 154{,}680{,}834 \\
        industrial\_07 & & 1{,}032{,}249      & 1{,}847{,}445      & 158{,}270{,}868 \\
        industrial\_08 & & 8{,}282{,}363      & 15{,}019{,}040     & 50{,}014{,}259 \\
        industrial\_09 & & 8{,}274{,}338      & 15{,}137{,}315     & 50{,}202{,}651 \\
        industrial\_10 & & 8{,}300{,}620      & 15{,}175{,}781     & 50{,}353{,}085 \\
        industrial\_11 & & 156{,}071{,}088    & 67{,}852{,}370     & 378{,}757{,}180 \\
        industrial\_12 & & 47{,}214{,}138     & 21{,}912{,}905     & 120{,}056{,}229 \\
        industrial\_13 & & 75{,}162{,}219     & 56{,}461{,}136     & 210{,}912{,}914 \\
        industrial\_14 & & 53{,}218{,}831     & 23{,}205{,}747     & 132{,}082{,}546 \\
        \hline
        \multicolumn{5}{l}{\footnotesize $^{\mathrm{a}}$ available from \url{https://zenodo.org/records/18311795} } \\
        \multicolumn{5}{l}{\footnotesize $^{\mathrm{b}}$ available from \url{https://zenodo.org/records/19205518}} \\
        \multicolumn{5}{l}{\footnotesize $^{\mathrm{c}}$ available from \url{https://openenergybenchmark.org/}}
         
    \end{tabular}

    \label{tab:instances}
\end{table}

\section{LP Solution Paradigms}
\begin{table}[h]
    \centering
    \footnotesize
    \caption{Comparison of shared-memory, distributed-memory, and GPU-accelerated solution paradigms for large-scale linear energy system optimization}
    \begin{tabular}{p{3cm}|p{4cm}|p{4cm}|p{4cm}}
 & \textbf{Shared-memory IPM} & \textbf{Distributed-memory IPM} & \textbf{GPU-accelerated FOM} \\ \hline
\textbf{Hardware\linebreak Requirement} 
& High-end workstation / server 
& Multi-node HPC cluster 
& Modern GPU hardware \\ \hline
\textbf{Memory Limit} 
& Bounded by single-node RAM 
& Aggregated cluster RAM 
& Bounded by VRAM capacity \\ \hline
\textbf{Solution Accuracy} 
& High / machine precision or user-defined tolerance 
& High / machine precision or user-defined tolerance 
& Medium / user-defined tolerance \\ \hline
\textbf{Parallel Scalability} 
& Limited 
& High 
& Massive \\ \hline
\textbf{Implementation\linebreak Effort} 
& Low; through standard interfaces 
& High; specialized solution needed 
& Moderate; through recent API support \\ \hline
\textbf{Structural\linebreak Dependency} 
& None (black-box) 
& High 
& None (matrix-vector based) \\ \hline
\textbf{Iteration\linebreak Convergence} 
& Fast quadratic 
& Quadratic; requires inter-node communication 
& Slow; first-order gradient steps \\ \hline
\textbf{Handling of Ill-conditioning} 
& Strong; handled via direct factori ation 
& Moderate; sensitive to partition quality 
& Robust; prevents crashes but slows convergence \\ \hline
\textbf{Cost of Entry} 
& Standard licensing / hardware 
& High (HPC infrastructure) 
& Medium (high-end GPUs) \\
& &
\multicolumn{2}{c}{access in academic computing environments }
\end{tabular}

    
    \label{tab:methods}
\end{table}
In the context of large-scale linear energy system optimization, standard shared-memory IPM represent the established state-of-the-art. Their prominence is rooted in a combination of high numerical reliability, the ability to achieve machine precision, and a mature ecosystem of standard interfaces that require minimal implementation effort from the modeller. By utilizing second-order Newton steps, these solvers exhibit fast quadratic convergence, typically reaching an optimal solution in a low, predictable number of iterations. Furthermore, they can effectively utilize multi-core processors through multithreading, but see diminishing performance gains between 10 and 20 threads, depending on the specific problem. Thus, the speedup through multithreading usually levels off at a sweet spot, which is much lower than the core counts (200+) on modern high-end server CPUs. However, they are still a highly performant and most accessible choice for a wide range of models that can fit within the memory constraints of a single high-end workstation or server node. However, in light of increasingly complex models with finer resolutions and the need to handle uncertainties, the memory limits of traditional IPMs on CPUs necessitates the consideration of emerging solution strategies, presented by distributed-memory IPMs and GPU accelerated FOMs. The main properties of the three solution paradigms are summarized in Table \ref{tab:methods}. 

Distributed-memory IPMs extend the reach of second-order methods by aggregating the RAM of multiple nodes in an HPC cluster. They aim to solve LPs by exploiting some underlying problem structure and parallelizing the sparse linear algebra of partial problems through distributed direct solvers. For these solvers in particular, the identification of a suitable structure is a crucial aspect, since the structure determines the required communication overhead between tasks, which can diminish the achievable parallelism. Here, PIPS-IPM++ \parencite{Rehfeldtetal22} is an example of a solver designed to exploit the inherent structure of linear energy system optimization problems, distribute blocks across multiple cores and nodes, and handle linking variables and constraints efficiently. While this approach offers high parallel scalability, it introduces a significant implementation burden, as it either requires the modeler to explicitly annotate the problem to manually define a block-angular structure or rely on automatic detection approaches, which can widely differ in quality for different problems. Furthermore, the reliance on inter-node communication and the sensitivity of the distributed system to partition quality can lead to variable numerical robustness compared to the direct factorization used in shared-memory environments.

Alternatively, GPU-accelerated FOMs have emerged as a paradigm shift for ultra-large-scale optimization. Unlike the second-order IPMs, FOMs rely on matrix-vector products that are ideally suited for the massive parallelism of modern GPU hardware \parencite{Applegate25, ZhangBoyd25, HPR-LP2025}. They can often handle problems with millions of variables.  FOMs are numerically robust and rarely crash due to ill-conditioning. However, they often require thousands of iterations to reach a solution. Consequently, they are usually restricted to a medium level of solution accuracy. 

The adoption of these advanced paradigms is heavily influenced by the divergent computational environments of academia and industry. In academic research, the availability of high-performance computing (HPC) clusters, either through institutional ownership or national research infrastructure, lowers the barrier to experimenting with GPUs and distributed-memory architectures. Conversely, in a commercial context, the adoption of distributed or GPU-accelerated solvers is governed by a stricter cost-benefit analysis. While cloud computing provides on-demand access to massive hardware resources, the additional engineering overhead and the potential for non-deterministic solution times in the case of FOMs, that arise from variability in shared hardware, and runtime scheduling, or complex partitioning requirements (in the case of distributed IPMs) often favor the continued use of vertically scaled, single-node shared-memory solvers until they become physically untenable. 

Ultimately, while standard IPMs remain the baseline for most applications, both distributed-memory and GPU-based approaches offer specialized pathways for overcoming the physical limitations of single-node computation.

 \section{Experimental Setup}
To evaluate the performance and scalability of the different solution paradigms, we conduct a comprehensive benchmarking suite across shared-memory, distributed-memory, and GPU-accelerated architectures. To ensure the reproducibility and fairness of the benchmarks, all solvers are maintained in their default configurations where possible. Specifically, we use each solver's internal presolve routines. While preprocessing can significantly affect solution times, as demonstrated on the oemofB3\_int instances by Koch et al. (\citeyear{KochURB2025}), the extensive discussion of these sensitivities is beyond our scope; our goal is to evaluate the raw scalability of the underlying algorithms. All experiments are subject to a 24-hour time limit. To ensure results are comparable we set the solver's parameters to a uniform convergence tolerance of $10^{-6}$ for both primal and dual feasibility. 

To establish a performance ceiling for traditional architectures, we define a Virtual Best Shared-Memory (VB IPM) baseline. This represents the fastest execution time achieved by any of the selected commercial solvers (COPT 8.0.3, CPLEX 22.1.2.0, Gurobi 13.0.0, FICO Xpress 9.8.0, and MOSEK 11.1.5) for a given instance.  This baseline represents the ``industry standard'' performance. We compare this against HiGHS 1.12.0 \parencite{Qi2018}, an open-source solver frequently used in the energy community. Tests for shared-memory IPM are conducted using 16 and 32 threads on single-node hardware\footnote{Intel(R) Xeon(R) Gold 6342 CPU} to assess how these models scale as more CPU cores are added. The crossover phase is deactivated. In energy modeling, crossover is typically used to move from a ``fuzzy'' interior solution to a ``crisp'' vertex solution; by deactivating it, we focus strictly on the computational effort required to reach the numerical optimum (see discussion of the effect in Section \ref{section:results}).

We utilize PIPS-IPM++ for distributed-memory multi-node experiments. This solver requires the model to be partitioned into a ``block-angular'' structure, which in energy systems typically reflects independent time-slices or regions linked by shared constraints (like storage levels or annual emission targets). As the names of variables and equations are not disclosed for most instances in the test set, partitioning into blocks must be based solely on the matrix structure. To this end, we use a hypergraph partitioning algorithm based on the Karlsruhe Hypergraph Partitioner (KaHyPar) \cite{gottesburen2024scalable} to identify equally sized groups of variables and equations, which are sparsely linked to other clusters. As a target for the partitioning, we utilize 80~blocks to allow a good fit of the number of blocks on the available HPC infrastructure. While problem-specific tuning of the block partitioning can further accelerate models \parencite{WETZEL2025101893}), we use a standardized default to maintain a fair comparison with other automated solvers. The experiments are conducted on two High-Performance Computing (HPC) systems: terrabyte (Intel-based, 2x40-core Intel Xeon Platinum 8380 processors with 1TB RAM per node) and DLR CARA (AMD-based, 2x64-core AMD EPYC 7702 processors with 1TB RAM per node). The experiments on CARA employed Panua Pardiso \cite{schenk2004solving} as the sparse linear solver, the experiments on terrabyte HSL-MA57 \parencite{Duff2004}.

We evaluate the performance of acceleration methods employing modern GPU hardware, specifically the NVIDIA H200, to assess its viability for ultra-large-scale energy systems. We define a Virtual Best FOM (VB FOM) using first-order methods (PDHG in Gurobi and Xpress, and PDLP in COPT) and compare the VB FOM baseline against the open source solution NVIDIA cuOpt 26.4.00a78\footnote{Due to a bug in the PSLP presolver, the released version 26.02 was not suitable for this study. We worked with NVIDIA and Daniel Cederberg (PSLP) to resolve this quickly and use this nightly build version.}.

 \section{Computational Results} \label{section:results}
\begin{table}[ht]
    \centering
    \caption{Solution Summary: Instances out of 46 solved to optimality/t-out/other status or not run} 
    \begin{tabular}{rrrrr}
        \toprule
          \textbf{VB IPM} & \textbf{HiGHS} & \textbf{PIPS} & \textbf{VB FOM} & \textbf{cuOpt}  \\
        \midrule
         43/ 3/ 0 & 8/ 22/ 16 & 14/ 0/ 32 & 42/ 4/ 0 & 40/ 1/ 5 \\
         \bottomrule
    \end{tabular}
    
    \label{tab:solution_summary}
\end{table}

Table \ref{tab:solution_summary} details the solvability of the test set within the prescribed time limits and precision thresholds across the evaluated solution methods. The VB IPM baseline finds optimal solutions for all but three instances. Thus, at least one of the commercial IPM solvers computes the optimal solution. Note that no single solver dominated across this diverse test set. Thus, the instance-wise virtual best approach yields significantly better overall results than any individual solver used.  In contrast, HiGHS achieves optimality for only 8 instances, hitting the time limit in 21 of them. Similarly, over the full test set PIPS-IPM++ demonstrates limited performance, solving only 14 of the benchmark instances (the best performance of PIPS-IPM++ on its two hardware/sparse solver variants for each instance). Regarding FOM solvers, both the commercial VB FOM and cuOpt consistently produce solutions, with only minor exceptions. Remarkably, none of the IPM solves industrial\_11 within the requested time limit, while cuOPT computes an optimal solution (see Table \ref{tab:all_times}). 

\begin{figure}[t]

    \centering
    \begin{subfigure}[b]{0.48\linewidth}
        \centering
    \includegraphics[width=\linewidth]{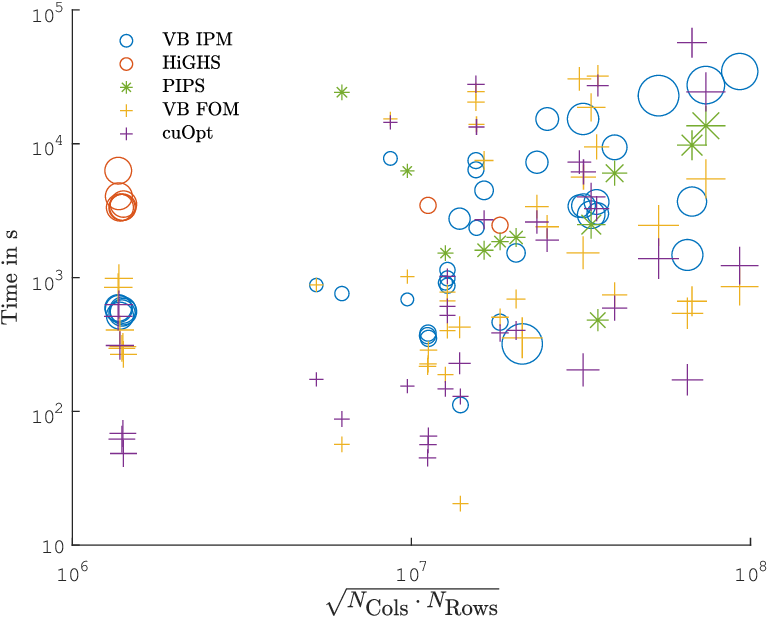}
    \subcaption{Solution times w.r.t. to constraint matrix properties. The total size is represented by $\sqrt{\text{N}_{\text{Cols}} \cdot \text{N}_{\text{Rows}}}$.} 
    \label{fig:solution_times}
    \end{subfigure}
    \hfill 
    \begin{subfigure}[b]{0.48\linewidth}
        \centering
         \includegraphics[width=\linewidth]{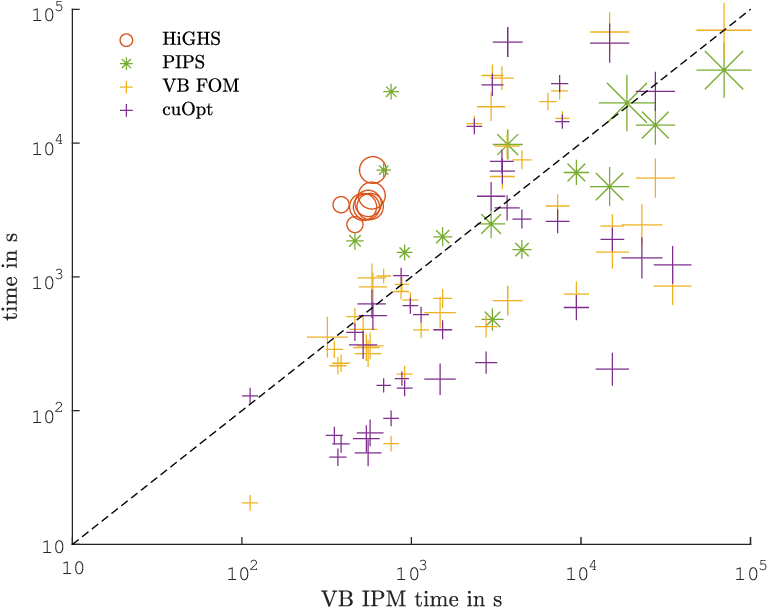}
         \subcaption{ Solution times w.r.t. to VB IPM, where the diagonal line ($y=x$) represents parity with VP IPM.}
         \label{fig:speedup}
         \end{subfigure}
         \caption{Summary on times to optimality. Results are plotted on a log-log scale and the size of the markers is scaled according to the number of non-zeros.}
         \label{fig:times}
    
\end{figure}

Despite these successes, several overarching computational challenges were observed. Many ESOM instances triggered numerical warnings and related premature barrier solver termination. While some solvers offer parameters for specialized support, permanent resolution often requires model-level reformulations and proper scaling of input data. Since these constraints are typically handled by the modeling framework, implementing such fixes requires expert knowledge of optimization techniques to override standard abstractions. For energy system modelers, this requires critical reflection on the choice of units for input data and the trade-off between numerical stability of the optimization problem and the necessity of optional constraints for answering the specific research question at hand. For the largest instance (BEAM\_4032\_11\_8\_CLI), several solvers failed to initiate the optimization process. The high number of non-zero elements led to memory exhaustion or internal integer overflows during the initial file-reading phase, highlighting a bottleneck in model instance ingestion for large-scale ESOMs. Detailed solution times or exit status if not optimal for each instance are provided in Table~\ref{tab:all_times} and aggregated in Fig.~\ref{fig:times}. Fig.~ \ref{fig:solution_times} relates solution times to overall instance size. Although instance size and time to solution appear to be related, there is a strong variability. Likewise, when comparing solution time to sparsity metrics, no clear relationship emerges. This suggests that the computational complexity is driven by the specific numerical conditioning of the constraints rather than simple matrix density. Comparing solution times compared to the VB IPM solution time, Fig.~\ref{fig:speedup} reveals that GPU-accelerated FOMs (VB FOM and cuOpt) offer transformative speedups for some of the largest instances. Interestingly, for the most challenging instances -- those requiring over $10^4$ seconds for the VB IPM -- the FOMs frequently achieve speedups of one to two orders of magnitude. However, it is important to contextualize these gains within the medium accuracy profile of these solvers. 

For example, we evaluated instance 288\_22\_8 on an NVIDIA Spark architecture, highlighting critical factors affecting solution quality. Initial results using the VB IPM solver reached a primal-dual interior solution in 1,820 seconds on this hardware. However, independent numerical validation of this point revealed primal infeasibilities of up to approximately $10^{-1}$ and substantial violations of the complementary slackness conditions. Achieving a high-precision solution required a transition to a crossover phase, which increased the total runtime to 10,417 seconds. After the crossover phase, primal infeasibility was reduced to $10^{-6}$. In contrast, testing with cuOpt demonstrated that simply tightening feasibility and optimality tolerances (from $10^{-4}$, over $10^{-6}$ to $10^{-8}$) led to a massive escalation in iterations and runtime (from 555, over 4050 to  9,727 seconds). Yet, the resulting improvements in primal and dual violations were marginal. Furthermore, the crossover phase for the GPU-based approach failed to converge within the designated time limits. 

However, not only the qualitative, but also the quantitative aspect is difficult to assess. Even on the same hardware, we notice a large performance variability. For example, when solving 288\_22\_8 with five different random seeds, we observe fluctuations of the total solving time of up to 48\,\% around the average time. Moreover, the memory consumption varies significantly: It happens that one run terminates normally with 20\,\% memory utilization, while another runs out of memory. We can nevertheless confirm that FOMs are less (main-)memory-demanding. In comparison to VB IPM, the memory footprint of cuOpt is reduced by a factor of 2.1 on average. This is an important practical aspect, as our larger instances require several hundred gigabytes of RAM. 

For the distributed memory solver PIPS-IPM++, a large number of solver errors could be traced back to internal integer overflows associated with a large number of linking elements between blocks, hinting at either an insufficient quality of the block-detection algorithm or model-specific formulations preventing the partitioning. Notably, all of the instances for which PIPS-IPM++ demonstrates a competitive performance are from models with a high temporal resolution such as REMix, oemof and PyPSA. The oemofB3\_int instances, which represent two-stage stochastic optimization problems, demonstrate the performance advantage that can be achieved when problem sizes increase, and the number of linking elements remains small. Going forward, the difference between partitioning methods that build on domain-specific knowledge and automatic annotation methods must be compared to identify heuristics for detecting block structures that are more suitable for PIPS-IPM++. Furthermore, the strong dependence of performance on the chosen direct solver highlights the need to integrate different solvers to better understand the performance bottlenecks associated with the various parts of the Schur complement decomposition method.

\section{Conclusion} \label{section:conclusion}
Our benchmarking of 46 large-scale ESOM instances demonstrates that no single solution paradigm dominates; performance depends heavily on the model instance's structural coupling. Shared-memory IPMs remain the baseline for high-accuracy solutions, typically reaching optimality within 24 hours. However, for block-angular structures -- common in two-stage stochastic or high-temporal-resolution models -- distributed-memory IPMs like PIPS-IPM++ leverage HPC parallelism to scale across hundreds of cores by exploiting limited linking variables.

Meanwhile, GPU-accelerated FOMs offer a paradigm shift for ultra-large-scale LPs where medium numerical accuracy is acceptable. In many energy applications, uncertainties in model assumptions and data projections exceed the solver’s numerical residuals. Consequently, extreme precision may be less critical for robust interpretation than the ability to compute solutions at previously unreachable scales.

Crucially, results suggest that computational complexity is driven by numerical conditioning and specific constraint coupling, such as investment decisions, rather than strictly size or matrix density.

Beyond algorithms, hardware requirements are substantial. IPMs are primarily RAM-constrained due to matrix factorization overhead, whereas FOMs shift the burden to high-end GPU performance. Accessing specialized hardware like NVIDIA H200 systems via cloud services remains expensive and limited, often requiring prior reservations. These findings reflect a current snapshot; however, both optimization algorithms and hardware architectures continue to evolve rapidly.

\section*{Acknowledgments}

The work for this article has been conducted in the Research Campus MODAL funded by the German Federal Ministry of Research, Technology and Space (BMFTR) (fund numbers 05M2025) and within the project PEREGRINE funded by the German Federal Ministry for Economic Affairs and Energy (BMWE) under grant number 03EI1082A. The authors gratefully acknowledge the scientific support and HPC resources provided by the German Aerospace Center (DLR). The HPC system CARA is partially funded by “Saxon State Ministry for Economic Affairs, Labour and Transport” and “Federal Ministry for Economic Affairs and Energy”. The authors gratefully acknowledge the computational and data resources provided through the joint high-performance data analytics (HPDA) project “terrabyte” of the German Aerospace Center (DLR) and the Leibniz Supercomputing Center (LRZ).

\section*{Appendix}

\setcounter{table}{0}
\renewcommand{\thetable}{A.\arabic{table}}

\begin{table}[H]
\centering
\caption{Summary of solution times and status in seconds. \\ \footnotesize cuOpt times marked * are run on an H100 machine and scaled to an H200 architecture by a factor of 0.873 derived from the instances available on both machines.}
\scalebox{0.94}{%
\footnotesize
\begin{filecontents*}{results/solution_times_table2.csv}
Name ; total_min_16thr ; total_min_32thr ; min_highs_hipo_16thr ; min_highs_hipo_32thr ; min_pips;total_min_fom ; min_cuopt_pdg
288\_22\_8                ;2755;3327; t-out           ; t-out           ; else            ;425;229
672\_22\_8\_cli           ;16918;15322; t-out           ; t-out           ; else            ;1530;204
BEAM\_2016\_11\_8         ;26949;22888;else;else; else            ;2451;1385
BEAM\_4032\_11\_8\_CLI    ; t-out           ; t-out           ;else;else; else            ; t-out           ; else           
ELMOD\_876\_10\_noVEnames ;19253;18657; t-out           ; t-out           ;19978; t-out           ; else           
ethos\_fine\_europe\_60tp-175-720ts ;10278;7769; t-out           ; t-out           ; else            ;15327; 14470*
genx-elec\_co2-15-168h    ;869;931; else            ; else            ; else            ;778;1021
genx-elec\_trex\_co2-15-168h ;1177;1139; else            ; else            ; else            ;401;522
genx-elec\_trex-15-168h   ;1065;986; t-out           ; t-out           ; else            ;670;609
industrial\_01            ;1542;1528; t-out           ; t-out           ;1993;690;403
industrial\_02            ;637;543;3456;3337; else            ;296;62
industrial\_03            ;685;571;3666;3358; else            ;305;68
industrial\_04            ;677;556;3777;3540; else            ;267;48
industrial\_05            ;686;584;4290;4051; else            ;984;628
industrial\_06            ;725;592;6573;6294; else            ;843;513
industrial\_07            ;586;518;3713;3338; else            ;406;311
industrial\_08            ;404;368; else            ; else            ; else            ;217;45
industrial\_09            ;384;400;3603;3462; else            ;226;56
industrial\_10            ;351;424; else            ; else            ; else            ;287;65
industrial\_11            ; t-out           ; t-out           ; t-out           ; t-out           ; --              ; t-out           ;59558
industrial\_12            ;3989;3440; t-out           ; t-out           ; else            ;5626; 6162*
industrial\_13            ;1497;1475; else            ; else            ; --              ;540;172
industrial\_14            ;4451;3671; else            ; else            ; else            ;9489; 3278*
JRC-EU-TIMES-dispatch-30-2016ts ;319;571; infeas.         ; infeas.         ; else            ;355; infeas.        
nt2030\_2030\_CY2009\_st1 ;464;482;2460;2551;1857;505;385
oemofB3\_int\_125         ;14738;21802; t-out           ; t-out           ;4728;67673; 55883*
oemofB3\_int\_250         ;69791;72362; else            ; else            ;35072;69738; else           
oemofB3\_int\_37          ;3005;3061; t-out           ; t-out           ;481;32004; 27228*
OTAI\_b-ts5\_cplex        ;3132;2950; else            ; else            ;2486;18695;4006
pypsa-de-elec-50-1h       ;952;911; t-out           ; t-out           ;1522;188;147
pypsa-de-sec-20-1h        ;3714;3415; t-out           ; t-out           ; else            ;30592; 7288*
pypsa-eur-elec-100-3h     ;747;686; t-out           ; t-out           ;6272;1015;154
remix\_nagsys\_cwe\_3h    ;8003;7274; t-out           ; t-out           ; else            ;3386;2601
remix\_nagsys\_cwe\_ineq\_1h ;29990;27453; t-out           ; t-out           ;13637;5463; 24392*
remix\_nagsys\_eu\_8h     ;17808;15310; t-out           ; t-out           ; else            ;2398;1906
remix\_nagsys\_eu\_ineq\_1h ; t-out           ; t-out           ; else            ; else            ; else            ; t-out           ; else           
remix\_nagsys\_fr\_1h     ;4742;4484; t-out           ; t-out           ;1600;7497; 2704*
remix\_unseen\_adalpert\_1h ;36476;34670; t-out           ; t-out           ; else            ;854;1229
remix\_unseen\_gunilpert\_1h ;11371;9399; t-out           ; t-out           ;6029;743;590
remix\_yssp\_disp\_488r\_1h ;3917;3698; t-out           ; t-out           ;9759;665; 56941*
spineopt-multiyear-invest-eu-case-study ;5078;2352; t-out           ; else            ; else            ;13980; 13343*
SWITCH-China-open-model-32-433ts ;141;112; else            ; else            ; else            ;20; 129*
TIMES-GEO-global-base-31-20ts ;816;759; else            ; else            ;24237;57; 88*
TIMES-GEO-global-netzero-31-20ts ;968;879; else            ; else            ; else            ;880; 173*
zen-garden-eur-PI-annual-emission-limit-28-100ts ;7757;6386; t-out           ; t-out           ; else            ;20443; t-out          
zen-garden-eur-PI-constrained-expansion-28-100ts ;9873;7485; else            ; else            ; else            ;24563;27814
\end{filecontents*}
\pgfplotstabletypeset[
    col sep=semicolon,
    trim cells=true,
    header=true,
    string type, 
    columns={[index]0,[index]1,[index]2,[index]3,[index]4,[index]5,[index]6,[index]7},
    display columns/0/.style={
        column name={Instance Name},
        column type=p{5.75cm},
        string type,
        postproc cell content/.append style={
        /pgfplots/table/@cell content/.add={\scriptsize}{},
        }
    },
    display columns/1/.style={column name={VB-IPM-16}, column type={>{\raggedleft\arraybackslash}p{1.2cm}}},
    display columns/2/.style={column name={VB-IPM-32}, column type={>{\raggedleft\arraybackslash}p{1.2cm}}},
    display columns/3/.style={column name={HiGHS-16}, column type={>{\raggedleft\arraybackslash}p{1.2cm}}},
    display columns/4/.style={column name={HiGHS-32}, column type={>{\raggedleft\arraybackslash}p{1.2cm}}},
    display columns/5/.style={column name={PIPS-IPM++}, column type={>{\raggedleft\arraybackslash}p{1.2cm}}},
    display columns/6/.style={column name={VB FOM}, column type={>{\raggedleft\arraybackslash}p{1.2cm}}},
    display columns/7/.style={column name={cuOpt}, column type={>{\raggedleft\arraybackslash}p{1.2cm}}},
    every head row/.style={before row=\toprule, after row=\midrule},
    every last row/.style={after row=\bottomrule},
]{results/solution_times_table2.csv}
\label{tab:all_times}
}
\end{table}
\section*{References}
\printbibliography[heading=none]

\end{document}